\documentstyle[amsfonts]{article}
\newtheorem{theorem}{Theorem}
\newtheorem{corollary}{Corollary}
\newtheorem{lemma}{Lemma}
\newtheorem{question}{Question}
\newcommand{\Pf}{{\cal P}(\omega)/fin}
\newcommand{\qed}{$\Box$}

\begin{document}
\title{Cardinal $\bf p$ and a theorem of Pelczynski}
\author{M. V. Matveev\\
{\small Department of Mathematics, University of California, Davis}\\
{\small Davis CA 95616, USA (address valid till June 30, 2000}\\
{\small e-mail: misha$\underline{\mbox{ }}$matveev@hotmail.com
}}
\date{}
\maketitle

\noindent
{\small{\bf Abstract.}
Are two compactifications of $\omega$ homeomorphic if their remainders
are homeomorphic?
For metrizable compactifications the question was answered
affirmatively by Pelzcynski.
Can the same happen for some non-metrizable remainders?
We consider the case when the remainder is $D^\tau$
for some uncountable $\tau$.
We show that the answer is affirmative if $\tau<\bf p$
and negative if $\tau=\bf c$.
We prove that every isomorphism between two subalgebras of $\Pf$
is generated by a permutation of $\omega$ provided these
subalgebras have independent basis of cardinality fewer than $\bf p$.
Also we consider some special dense countable subsets in $D^\tau$.
}

\medskip
\noindent
{\bf AMS subject Classification:} Primary 54D35; Secondary 54D65,
54A20, 03E10, 06E05

\medskip
\noindent
{\bf Keywords:} compactification, remainder, Boolean algebra,
independent, cardinal $\bf p$, cardinal $\bf b$,
sequentially separable.

\section{Introduction}

In 1965 Pelczynski proved the following:

\begin{theorem}\label{ofPel}
{\em \cite{Pel}}
Let $b_1\omega$ and $b_2\omega$ be two metrizable compactifications
of $\omega$.
If $b_1\omega\setminus\omega$ and $b_2\omega\setminus\omega$ 
are homeomorphic, then $b_1\omega$ and $b_2\omega$ are homeomorphic.

Moreover, every continuous mapping 
$f:b_1\omega\setminus\omega\to b_2\omega\setminus\omega$
can be extended to a continuous mapping $F:b_1\omega\to b_2\omega$
so that $F(\omega)\subset\omega$ and 
if $f$ is surjective, then so is $F$.
\end{theorem}

(See also Terasawa's paper \cite{Terasawa} for an alternative proof 
and related results).

In this paper we are interested in the question whether
Pelczynski's theorem can be extended to some non-metrizable compacta.
We obtain the following particular results:

\begin{theorem}\label{positive}
If $\tau<\bf p$, then all compactifications of $\omega$
with the remaider homeomorphic to $D^\tau$ are homeomorphic 
to each other.
\end{theorem}

\begin{theorem}\label{negative}
There are $2^{\bf c}$ pairwise non-homeomorphic compactifications
of $\omega$ with the remainder homeomorphic to $D^{\bf c}$.
\end{theorem}

(The reader is refered to \cite{vD} or \cite{Vau}
for the definition of cardinal $\bf p$
and other small cardinals and $^*$-terminology).
The following natural questions remain open.

\begin{question}{\em
Can one construct within ZFC a non-metrizable compact space $X$
such that all compactifications of $\omega$ with the remainder
homeomorphic to $X$ are homeomorphic to each other?
}
\end{question}

\begin{question}\label{two}
{\em
For which cardinal numbers $k$, $1<k<2^{{\bf c}}$, is there a compact $X$
such that there are exactly $k$ many pairwise non-homeomorphic
compactifications of $\omega$ with the remainder homeomorphic to $X$?
}
\end{question}

Restriction $<2^{{\bf c}}$ in Question~\ref{two} is natural:
one cannot construct more than $2^{{\bf c}}$ pairwise non-homeomorphic
compactifications of $\omega$.

In the  proof of Theorem~\ref{positive} we essentially use the fact 
that the Boolean algebra of clopen subsets of $D^\tau$
has independent basis. So it is not clear how to extend Theorem~\ref{positive}
to other compacta of small weight.

The basic idea of the proof of Theorem~\ref{negative}
is the notion of sequential separability (or, rather, its various negations)
which we discuss in Section~\ref{SeqSep}.

An interesting question was rised by R. Williams in the review of 
Pelczynski's paper (MR 32 $\#$4659): can the assignment $f\to F$
in Theorem~1 be made functorial (i.e. such that $FG$
always corresponds to $fg$)?

\section{Proof of Theorem 2}

$\Pf$ denotes the Boolean algebra of classes of subsets of $\omega$
modulo the ideal of finite sets.
A family $\cal A$ of subsets of $\omega$
is independent provided for every
$n,m\in\omega$ and every distinct $A_1,\dots, A_{n+m}\in \cal A$ the
intersection
$A_1\cap\dots\cap A_n\cap(\omega\setminus A_{n+1})\cap(\omega\setminus
A_{n+m})$ is nonempty.
In the same sense we speak about the independence of subfamilies of $\Pf$.
We say that an isomorphism $f: {\cal A}\to{\cal B}$ between two
subalgebras
${\cal A},{\cal B}\subset\Pf$ is generated by a permutation $\pi$ of
$\omega$ provided
for every $a\in\cal A$ and every $A\in a$, $\pi(A)\in f(a)$.

It is clear that any subalgebras of $\Pf$ having independent basis of
the same cardinality are isomprphic.

\begin{lemma}\label{BAInd}
Let $f:{\cal A}\to\cal B$ be an isomorphism between subalgebras
${\cal A},{\cal B}\subset\Pf$.
If $|{\cal A}|<\bf p$ and $\cal A$ has an independent basis
then $f$ can be generated by a permutation of $\omega$.
\end{lemma}

{\bf Proof:}
Denote $\tau=|{\cal A}|$. We can assume that $\tau\geq\omega$.
Let ${\cal A}_0$ be an independent basis for $\cal A$.
Then $|{\cal A}_0|=\tau$ and we can write  
${\cal A}_0=\{A_\alpha:\alpha<\tau\}$.
It is clear that $\{f(A_\alpha):\alpha<\tau\}$ is an
independent basis for $\cal B$.
For each $\alpha<\tau$, choose $a_\alpha\in A_\alpha$ and $b_\alpha\in
f(A_\alpha)$.
For $n\in\omega$ and $\alpha<\tau$, put
$$
S_{n,\alpha}=\left\{\begin{array}{l}
b_\alpha,\mbox{ if }n\in a_\alpha\\
\omega\setminus b_\alpha\mbox{, otherwise}.
\end{array}\right.
$$
Then the family
${\cal S}_n=\{S_{n,\alpha}:\alpha<\tau\}$
has sfip and, since $\tau<\bf p$, a pseudointersection, say $H_n$.

Next, for each $\alpha<\tau$ we define a function $h_\alpha:\omega\to
\omega$ by
$h_\alpha(n)=\min\{m: H_n\setminus\{0,\dots,m-1\}\subset S_{n,\alpha}\}$.
Since $\tau<{\bf p}\leq{\bf b}$, there is a strictly increasing function
$g_0:\omega\to\omega$ such that $g_0\geq^* h_\alpha$ for all $\alpha<\tau$.
Define inductively a strictly increasing function $\pi_0:\omega\to\omega$
such that $\pi_0\geq g_0$ and $\pi_0(n)\in H_n$ for all $n$.
Then $\pi_0$ is an injection of $\omega$ to $\omega$ such that
$\pi_0(a_\alpha)\subseteq^* b_\alpha$ 
and $\pi_0(\omega\setminus a_\alpha)\subseteq^* \omega\setminus b_\alpha$ 
for all $\alpha<\tau$.
Working simetrically, we obtain an injection $\pi_1:\omega\to\omega$
such that $\pi_1(b_\alpha)\subseteq^* a_\alpha$ 
and $\pi_1(\omega\setminus b_\alpha)\subseteq^* \omega\setminus a_\alpha$ 
for all $\alpha\in\tau$. 
By the usual Cantor-Bernstein argument 
(see, for example \cite{KM}, Theorem~5.5.2) there are decompositions
$\omega=M_1\cup M_2$ and $\omega= N_1\cup N_2$ such that
$\pi_0(M_1)=N_1$ and $\pi_1(N_2)=M_2$.
Then the function $\pi$ defined by the role
$$
\pi(n)=\left\{\begin{array}{l}
\pi_0(n),\mbox{ if }n\in M_1\\
\pi_1^{-1}(n),\mbox{ if }n\in M_2
\end{array}\right.
$$
is a permutation of $\omega$ and satisfies the condition
$\pi(a_\alpha)=^* b_\alpha$ for all $\alpha<\tau$.
So $f$ is generated by $\pi$.
\qed

\bigskip

For a space $X$, let $CO(X)$ denote the family of all clopen sets in $X$.
The next lemma follows directly from normality and compactness.

\begin{lemma}\label{omext}
Let $b\omega$ be a compactification of $\omega$ with a zero-dimensional
remainder $X=b\omega\setminus\omega$. 
Then there is a (unique up to $^*$)
mapping that
assigns to each $Y\in CO(X)$ a subset $g(Y)\subset\omega$
so that the set $Y'=Y\cup g(Y)$ is clopen in $b\omega$.
\end{lemma}

Now we prove Theorem~\ref{positive}.
Let $b_1\omega=(\omega\cup D^\tau,{\cal T}_1)$ and
$b_2\omega=(\omega\cup D^\tau,{\cal T}_2)$
be two compactifications of $\omega$ with remainder $D^\tau$.
Further, let $g_1,g_2:CO(D^\tau)\to {\cal P}(\omega)$
be mappings from Lemma~\ref{omext} corresponding to $b_1$ and $b_2$.
We will consider them as mappings to $\Pf$.
Note that ${\cal C}=\{\{f\in D^\tau:f(\alpha)=0\}:\alpha<\tau\}$
is an independent basis basis for the Boolean algebra $CO(D^\tau$. 
Then the subalgebras $\cal A$ and $\cal B$ of $\Pf$
generated by $g_1({\cal C})$ and $g_2({\cal C})$
also have independent basis and cardinality $\tau$.
It is clear that they are isomorphic.
By Lemma~\ref{BAInd} this isomorphism is generated
by some permutation, say $\pi$, of $\omega$.
Then the mapping $F:b_1\omega\to b_2\omega$
defined by the formula
$$
F(x)=\left\{\begin{array}{l}
x,\mbox{ if }x\in D^\tau\\
\pi(x),\mbox{ if }x\in\omega
\end{array}\right.
$$
is a homeomorphism.
\qed

\section{Sequential separability, its negations\\
and variations of cardinal $\bf p$}\label{SeqSep}

A space $X$ is {\em sequentially separable} provided there is 
a countable subspace $Y\subset X$ sucht that every point of $X$ is the limit
of a sequence from $Y$ ($Y$ is said to be {\em sequentially dense} in $X$).
Wilansky noted that the product of $\bf c$ sequentially separable
spaces need not be sequentially separable and asked if the product
of fewer than $\bf c$ sequentially separable spaces is
sequentially separable \cite{Wilansky}. 
Tall showed \cite{Tall} that 
the answer to this question depends on additional set-theoretic 
assumptions.
In fact, Tall proved that MA($k$) for $\sigma$-centered posets
implies that every product of fewer than $k$ sequentially separable
spaces is sequentially separable. 
Similar results were obtained in \cite{FK}.
Further, Bell proved (\cite{Bell}, see also \cite{Vau}, p.~201)
that ${\bf p}={\bf m}_{\sigma-\mbox{{\small centered}}}$ where
${\bf m}_{\sigma-\mbox{{\small centered}}}$ is the minimal cardinal $k$
such that ``MA($k$) for $\sigma$-centered posets'' fails.
So now Tall's theorem can be restated as follows.

\begin{theorem}
Every product of fewer than $\bf p$ sequentially separable spaces
is sequentially separable.
\end{theorem}

\begin{corollary}\label{trivial} 
$D^\tau$ is sequentially separable for all $\tau<\bf p$.
\end{corollary}

This particular case of Tall's theorem can be derived also from
the inequality ${\bf p}={\bf p}_{\bf p}={\bf p}_\chi$
(\cite{vD}, Theorem~6.2) where
${\bf p}_{\bf p}=\min\{k:D^k$  is not subsequential$\}$
and
${\bf p}_\chi=\min\{k:$ there exists a regular non-subsequential
space of character $k\}$
(subsequential, in van~Douwen's terminology, means that
if $x\in\overline{A}$ and $A$ is countable, then there is a sequence 
converging from $A$ to $x$).
So $D^k$ with $k<\bf p$ is better than just sequentially separable:
every dense countable subspace of $D^k$ is sequentially dense 
in it.\footnote{
By the way, ``countable'' in this statement is essential. 
Indeed, a $\Sigma$-product is dense but not sequentially dense
in $D^k$, moreover, it is sequentially closed.
}

Consider the following cardinal numbers.

\bigskip\noindent
${\bf p}_1=\min\{\tau:$ there is a dense countable $Y\subset D^\tau$
which is not sequentially 

\hskip 1cm dense$\}$,

\bigskip\noindent
${\bf p}_2=\min\{\tau: D^\tau$ is not sequentially separable$\}$,

\bigskip\noindent
${\bf p}_3=\min\{\tau:$ there is a dense countable $Y\subset D^\tau$
such that no nontrivial 

\hskip 1cm sequence from $Y$ converges in $D^\tau\}$. 

\bigskip
Then ${\bf p}\leq{\bf p}_1\leq{\bf p}_2\leq{\bf c}$
and ${\bf p}_1\leq{\bf p}_3\leq{\bf c}$.
The inequality ${\bf p}\leq{\bf p}_1$
follows from Corollary~\ref{trivial}; 
the inequality ${\bf p}_3\leq{\bf c}$
follows from Theorem~\ref{badsubspaces} below.
The rest is trivial.

Further, cardinals ${\bf p}_1$
and ${\bf p}_3$ can be expressed in a purely set-theoretic way.
Indeed, it is easy to see that

\bigskip\noindent
${\bf p}_1=\min\{\tau:$ there is an independent family
$\cal A$ of subsets of $\omega$ such 

\hskip 1cm that $|{\cal A}|=\tau$ and $\cal A$ does not have infinite
pseudointersection$\}$.

\bigskip
Further, let $\cal A$ and $\cal B$ be families of subsets of $\omega$.
We say that $\cal B$ is a {\it partial inversement} of $\cal A$
if there is a bijection $h:{\cal A}\to{\cal B}$
such that for each $A\in\cal A$ the set $h(A)$ equals either $A$ or
$\omega\setminus A$.
Then it is clear that 

\bigskip\noindent
${\bf p}_3=\min\{\tau$: there is an independent family $\cal A$ of subsets
of $\omega$ such that 

\hskip 1cm $|{\cal A}|=\tau$ and no partial inversement of $\cal A$
has nonempty 

\hskip 1cm pseudointersection$\}$.

\bigskip
Peter Nyikos noted \cite{Nyik} that in fact ${\bf p}_1=\bf p$.
Nyikos' proof goes as follows: let $\{I_\alpha:\alpha<{\bf p}\}$
be an independent family of subsets of $\omega$
and let $\{F_\alpha:\alpha<{\bf p}\}$ be a free filter on
$\omega$ with sfip but no infinite pseudointersection.
Then $\{I_\alpha\times F_\alpha:\alpha<{\bf p}\}$
is a ${\bf p}_1$-witnessing independent family of subsets
of $\omega\times\omega$. 

Rothberger proved in \cite{Rot} that $D^{\omega_1}$
is sequentially separable iff there is a Q-set of cardinality 
$\omega_1$.
In \cite{Lo} this result is extended as follows: $D^\tau$ is sequentially separable iff
there is a Q-set of cardinality $\tau$;
in other words, ${\bf p}_2=\bf q$ where $\bf q$ is the minimal
cardinal $\tau$ such that there is no Q-set of cardinality $\tau$
(cardinal $\bf q$ was introduced in \cite{GN}).
Since the existence of Q-sets of cardinality $\geq\bf p$
is consistent with ZFC (a model is constructed in \cite{FM} in which
there is a Q-set of cardinality $\omega_2$ and a non-Q-set of cardinality 
$\omega_1$
it is consistent also that ${\bf p}_1<{\bf p}_2$.
The question which of strict inequalities
${\bf p}_1<{\bf p}_3<{\bf c}$
can be consistently true remains open
as well as the relationship between ${\bf p}_2$ and ${\bf p}_3$
is not clear.

\begin{theorem}\label{badsubspaces}
There is a family ${\cal Y}=\{Y_\alpha:\alpha<2^{\bf c}\}$
of dense countable subspaces $Y_\alpha\subset D^{\bf c}$
such that: 

(1) no nontrivial sequence from $Y_\alpha$ converges in $D^{\bf c}$
and

(2) $Y_\alpha$ and $Y_{\alpha'}$ are disjoint and non-homeomorphic
as soon as $\alpha\neq\alpha'$.
\end{theorem}

(Trivial means eventially constant).

{\bf Proof:}
For convinience, we replace $\bf c$ with $\Bbb R$.
Two filters on $\omega$ have the same type if there is a permutation
of $\omega$ which transforms one of them into the other.
It is well known that there are $2^{\bf c}$ types of ultrafilters on 
$\omega$.
We construct $Y_\alpha$ by induction.
Suppose $\alpha<2^{\Bbb R}$ and $Y_\beta$ have been constructed for all
$\beta<\alpha$.
Denote $M_\alpha=\cup\{Y_\beta:\beta<\alpha\}$.
Then $|M_\alpha|<2^{\bf c}$.
Further, denote ${\cal F}_\alpha$ the set of all types of filters 
of neighbourhoods of $x$ in $Z$ where $x\in Z\subset Y_\beta$
for some $\beta<\alpha$
(if $\alpha=0$ then $M_\alpha=\emptyset$ and ${\cal F}_\alpha=\emptyset$).
Then $|{\cal F}_\alpha|<2^{\bf c}$ and thus we can pick
an ultrafilter $u_\alpha$ not in ${\cal F}_\alpha$.
Now we need a lemma.

\begin{lemma}
If $M\subset D^{\Bbb R}$ and $|M|<2^{\bf c}$,
then there is an embedding $i:\beta\omega\to D^{\Bbb R}$
such that $i(\beta\omega)\cap M=\emptyset$.
\end{lemma}

{\bf Proof of lemma:} We have $|pr_{{\Bbb R}^+}(M)|<2^{\bf c}=|D^{{\Bbb R}^+}|$
(we denote ${\Bbb R}^+=\{x\in{\Bbb R}: x\geq 0\}$ and
${\Bbb R}^-=\{x\in{\Bbb R}: x<0\}$).
Pick a point $z\in D^{{\Bbb R}^+}\setminus pr_{{\Bbb R}^+}M$.
Denote $Z=\{x\in D^{{\Bbb R}}:pr_{{\Bbb R}^+}(x)=z\}$.
Then $Z\cap M=\emptyset$ and $Z$ is homeomorphic to $D^{\bf c}$
and thus it contains a copy of every zero-dimensional space
of weight $\leq c$, in particular, a copy of $\beta\omega$.
\qed

Now we go on with the proof of the theorem.
By lemma, there is an embedding $i_\alpha:\beta\omega\to D^{\Bbb R}$
such that $i_\alpha(\beta\omega)\cap M_\alpha=\emptyset$.
Denote $m_\alpha=i_\alpha(u_\alpha)$.
Our $Y_\alpha$ will be of the form
$Y_\alpha=i_\alpha(\omega)\cup\{m_\alpha\}\cup Z_\alpha$
where $Z_\alpha$ is yet to be constructed.
This will be made in several steps.

First, consider the set $T$ of all such functions $f\in D^{\Bbb R}$ that $f$:

- takes each of the two values, $0$ and $1$, $\bf c$ times
on the set ${\Bbb R}^-$,

- has only finitely many points of discontinuity;

- all points of discontinuity are in $\Bbb Q$.

\noindent 
It is easy to see that $T$ is countable and dense in $D^{\Bbb R}$.

Next, enumerate all infinite, nowhere dense subsets of $T$ as
$\{N_\gamma:\gamma<{\bf c}\}$. 
For every $\gamma<{\bf c}$ fix a decomposition
$N_\gamma=N_\gamma^0\cup N_\gamma^1$ where both
$N_\gamma^0$ and $N_\gamma^1$ are infinite.
Also, enumerate the set ${\Bbb R}^+$ as
${\Bbb R}^+=\{q_\gamma:\gamma<{\bf c}\}$.
For every $t\in T$ we define $\tilde{t}\in D^{\Bbb R}$
as follows: for every $r\in\Bbb R$ put
$$
\tilde{t}(r)=\left\{\begin{array}{l}
0,\mbox{ if }r=q_\gamma\mbox{ and }t\in N_\gamma^0\\
1,\mbox{ if }r=q_\gamma\mbox{ and }t\in N_\gamma^1\\
t(r),\mbox{ otherwise.}
\end{array}\right.
$$
Put $\tilde{T}=\{\tilde{t}:t\in T\}$.
Then $\tilde{T}$ is countable.

We claim that $\tilde{T}$ is dense in $D^{\Bbb R}$.
Let $O=O_{r_1\dots r_n}^{i_1\dots i_n}=\{f\in D^{\Bbb R}:
f(r_1)=i_1,\dots,f(r_n)=i_n\}$ be a basic open set in $D^{\Bbb R}$.
Some of $r_1,\dots, r_n$, say 
$r_{k_1}=q_{\gamma_1},\dots, r_{k_s}=q_{\gamma_s}$ may be in 
${\Bbb R}^+$.
Note that the set $T\setminus(N_{\gamma_1}\cup\dots\cup N_{\gamma_s})$
is dense in $D^{\Bbb R}$.
Pick $t\in(T\setminus(N_{\gamma_1}\cup\dots\cup N_{\gamma_s}))\cap O$.
But then also $\tilde{t}\in O$ since $t$ and $\tilde{t}$ take the same
values at coordinates $q_{\gamma_1}$,..., $q_{\gamma_s}$.
So $\tilde{T}$ is dense in $D^{\Bbb R}$.

Now we claim that $\tilde{T}$ does not contain non-trivial sequences 
converging in $D^{\Bbb R}$.
Suppose that $\xi$ is a non-trivial sequence from $T$ converging 
in $D^{\Bbb R}$. 
Then the set $\tilde{S}=\{\xi(n):n\in\omega\}$
is infinite, and hence so is the set
$S=\{t\in T:\tilde{t}\in\tilde{S}\}$.
Further, $S$ contains an infinite nowhere dense set, say $N$.
But then $N=n_\gamma$ for some $\gamma<\bf c$,
and we have $\tilde{t}(\gamma)=0$ for all $t\in N_\gamma^0$ and
$\tilde{t}(\gamma)=1$ for all $t\in N_\gamma^1$.
This means that the sequence $\xi$ contains both infinitely many
elements taking value $0$ at coordinate $\gamma$ and
infinitely many elements taking value $1$ at coordinate $\gamma$.
Hence $\xi$ cannot converge.

Also we note that every point from $\tilde{T}$ takes each of the two
values, $0$ and $1$, $\bf c$ many times.
Indeed, $t$ and $\tilde{t}$ do not differ on the coordinates from
${\Bbb R}^-$.

Every permutation $\pi:{\Bbb R}\to\Bbb R$ induces
an autohomeomorphism $h_\pi:D^{\Bbb R}\to D^{\Bbb R}$
by the formula $(h_\pi(x))(r)=x(\pi^{-1}(r))$.
It is clear that $h_\pi(\tilde{T})$ has all properties of $\tilde{T}$
(is dense in $D^{\Bbb R}$, countable, and does not contain
nontrivial sequences converging in $D^{\Bbb R}$).
It remains to find such a permutation $\pi$ that $h_\pi(\tilde{T})$
does not intersect $M_\alpha$.
By induction it is easy to find for all $t\in T$ subsets 
$A_t, B_t\subset\Bbb R$ such that:

$|A_t|=|B_t|=\bf c$ for every $t\in T$,

$A_t\cap B_t=\emptyset$ for every $t\in T$,

$A_t\cap B_{t'}=\emptyset=A_t\cap A_{t'}$ whenever $t\neq t'$,

$\tilde{t}(a)=0$ for all $t\in T$ and $a\in A_t$,

$\tilde{t}(b)=1$ for all $t\in T$ and $b\in B_t$.

\noindent
Let $t\in T$. Denote ${\cal P}_t$
the set of all such permutations $\pi$ of $\Bbb R$
that $\pi(A_t)=B_t$, $\pi(B_t)=A_t$ and $\pi(a)=a$
for all $a\not\in A_t\cup B_t$.
It is clear that $|{\cal P}_t|=2^{\bf c}$.
Further, if $\pi,\pi'\in{\cal P}_t$ and $\pi\neq\pi'$, then
$pr_{A_t\cup B_t}(h_\pi(\tilde{t}))\neq pr_{A_t\cup B_t}(h_{\pi'}(\tilde{t}))$.
So $|\{pr_{A_t\cup B_t}(h_\pi(\tilde{t})):\pi\in{\cal P}_t\}|=2^{\bf c}$,
and since $|M_\alpha|<2^{\bf c}$
there is a $\pi_t\in{\cal P}_t$ such that 
$pr_{A_t\cup B_t}\left(h_{\pi_t}(\tilde{t})\right)\not\in 
pr_{A_t\cup B_t}(M_\alpha)$ ($*$).
The function $\pi:{\Bbb R}\to{\Bbb R}$ defined by formula
$$
\pi(a)=\left\{\begin{array}{l}
\pi_t(a),\mbox{ if }a\in A_t\cup B_t\mbox{ for some }t\in T\\
a,\mbox{ otherwise}
\end{array}\right.
$$
is a permutation.
It follows from ($*$) that 
$\pi(\tilde{T})\cap M_\alpha=\emptyset$.

Put $Z_\alpha=\pi(\tilde{T})$.
Then $Z_\alpha\cap M_\alpha=\emptyset$,
hence $Y_\alpha\cap M_\alpha=\emptyset$,
hence $Y_\alpha\cap Y_\beta=\emptyset$ for all $\beta<\alpha$.

The fact that $Y_\alpha$ is not homeomorphic to $Y_\beta$
for all $\beta<\alpha$ follows fro the inclusion
$Y_\alpha\supset i_\alpha(\omega)\cup \{m_\alpha\}$.

Last, if $Y_\alpha$ would contain a nontrivial sequence converging 
in $D^{\Bbb R}$, then this sequence would contain either 
a nontrivial subsequene from $Z_\alpha$, or one from
$i_\alpha(\omega)$.
The first is impossible by the construction of $Z_\alpha$,
the second - by the properties of \v{C}ech-stone compactification.
\qed

\section{Proof of Theorem 3 and more}

Let $Y\subset D^{\bf c}$ be one of $Y_\alpha$ from Theorem~\ref{badsubspaces}.
The following is a restriction of Tkachuk's construction from \cite{Tka}.
Fix a decomposition $\omega=\cup\{\omega_y:y\in Y\}$
where $\omega_y$ is infinite for every $y\in Y$.
Put $b_Y\tilde{\omega}=
(D^{\bf c}\times\{\omega\})\cup\tilde{\omega}\subset
D^{\bf c}\times(\omega+1)$
where $\tilde{\omega}=\cup\{\{y\}\times\omega_y:y\in Y\}$.
Then $b_Y\tilde{\omega}$ is compact; $\tilde{\omega}$ is countable,
dense in $b_Y\tilde{\omega}$ and consists of isolated points.
So $b_Y\tilde{\omega}$ can be considered as a compactification of $\omega$
with the remainder homeomorphic to $D^{\bf c}$.
Further, it is easy to see that
$\{x\in D^{\bf c}:\mbox{ there is a sequence converging from }
\tilde{\omega}\mbox{ to }x\}=Y$.
This fact and Theorem~\ref{badsubspaces} imply that 
$b_{Y_\alpha}\tilde{\omega}$ and $b_{Y_{\alpha'}}\tilde{\omega}$
are not homeomorphic whenever $\alpha\neq\alpha'$. \qed

\medskip
We conclude the paper with the following modification of Theorem~3.

\begin{theorem}
There are infinitely many pairwise non-homeomorphic
compactifications of $\omega$
with the remainder homeomorphic to $D^{{\bf p}_3}$.
\end{theorem}

{\bf Proof:}
Recall that the Alexandroff Duplicat $AD(X)$ of a topological space $X$
is the set $X\times\{0,1\}$ topologized by declaring the points of 
$X\times\{1\}$ to be isolated while basic neighbourhoods of points of 
$X\times\{0\}$ take the form $(U\times\{0,1\})\setminus\{(x,1)\}$
where $x\in X$ and $U$ is a neighbourhood of $x$ in $X$.
For $Y\subset X$ denote $AD_Y(X)=(X\times\{0\})\cup(Y\times\{1\})$
with the topology inhereted from $AD(X)$.
Let $Y\subset D^{{\bf p}_3}$ be dense, countable and
such that no non-trivial sequence
from $Y$ converges in $D^{{\bf p}_3}$.
Let $n\in\Bbb N$.
Choose $n$ distinct points $p_1,\dots,p_n\in D^{{\bf p}_3}$.
Further, let $A_1=N_1\cup\{a_1\}$,..., $A_n=N_n\cup\{a_n\}$
be $n$ disjoint copies of the convergent sequence
(the limit points are $a_1,\dots, a_n$).
Denote $b_{nY}\omega$ the quotient space of the discrete sum
$A_1\cup\dots\cup A_n\cup AD_Y D^{{\bf p}_3}$ obtained by
identifying $a_i$ with $(p_i,0)$ for all $1\leq i\leq n$.
Then $b_{nY}\omega$ can be considered as a compactification of 
$N_1\cup\dots\cup N_n\cup(Y\times\{1\})\sim\omega$
with the remainder $D^{{\bf p}_3}\times\{0\}\sim D^{{\bf p}_3}$.
It is clear that exactly $n$ points in $D^{{\bf p}_3}$ are limits
of sequences from $\omega$.
So compactifications $b_{nY}\omega$ are not homeomorphic
to each other for distinct $n$.
\qed

\bigskip
{\bf Asknowlegement.}
The paper was written when the author was visiting University of 
California, Davis. The author is grately indebted to colleagues
form UC Davis for their hospitality.

Also the author expresses his gratitude to Ronnie Levy for
numerous usefull discussions and to Peter Nyikos and Joseph Lo
for usefull discussions during Spring Topological Conference in March 2000.


\end{document}